\documentclass[10pt,fleqn,twoside]{article} %this is important to force latex!
\newcommand{\miniheader}{
  \documentclass[10pt,fleqn,twoside]{article}
  \usepackage{palatino}
  \usepackage{amsmath}
  \usepackage{amssymb}
  \usepackage{amsfonts}
  \usepackage{amsthm}
  \usepackage{eucal}
  \usepackage{graphicx}
  \usepackage{color}

  \usepackage{fancyhdr}
  \renewcommand{\headrulewidth}{.0pt}\renewcommand{\footrulewidth}{.0pt}\cfoot{}
  \fancyhead[OL,EC]{\it\theauthor---\today}
  \fancyhead[ER]{\leftmark}
  \fancyhead[OR,EL]{\thepage}
  \fancyfoot[EL,OR]{}

  \usepackage[round]{natbib}
  \bibliographystyle{abbrvnat}

  \graphicspath{{pics/}{figs/}{~/write/tex/pics/}{~/write/tex/figs/}{~/teaching/pics-all/}}
  \usepackage{geometry}
  \geometry{a4paper,hdivide={35mm,*,35mm},vdivide={35mm,*,35mm}}
  \renewcommand{\baselinestretch}{1.1}

  \renewcommand{\a}{\alpha}
  \renewcommand{\b}{\beta}
  \renewcommand{\d}{\delta}
    \newcommand{\D}{\Delta}
    \newcommand{\e}{\epsilon}
    \newcommand{\g}{\gamma}
    \newcommand{\G}{\Gamma}
  \renewcommand{\l}{\lambda}
  \renewcommand{\L}{\Lambda}
    \newcommand{\m}{\mu}
    \newcommand{\n}{\nu}
    \newcommand{\N}{\nabla}
  \renewcommand{\k}{\kappa}
  \renewcommand{\o}{\omega}
  \renewcommand{\O}{\Omega}
    \newcommand{\p}{\varphi}
  \renewcommand{\P}{\Phi}
  \renewcommand{\r}{\varrho}
    \newcommand{\s}{\sigma}
  \renewcommand{\S}{\Sigma}
  \renewcommand{\t}{\theta}
    \newcommand{\T}{\Theta}
    \newcommand{\x}{\xi}
    \newcommand{\X}{\Xi}
    \newcommand{\Y}{\Upsilon}
    \newcommand{\z}{\zeta}

  \renewcommand{\AA}{{\cal A}}
    \newcommand{\BB}{{\cal B}}
    \newcommand{\CC}{{\cal C}}
    \newcommand{\cc}{{\cal c}}
    \newcommand{\DD}{{\cal D}}
    \newcommand{\EE}{{\cal E}}
    \newcommand{\FF}{{\cal F}}
    \newcommand{\GG}{{\cal G}}
    \newcommand{\HH}{{\cal H}}
    \newcommand{\II}{{\cal I}}
    \newcommand{\KK}{{\cal K}}
    \newcommand{\LL}{{\cal L}}
    \newcommand{\MM}{{\cal M}}
    \newcommand{\NN}{{\cal N}}
    \newcommand{\oNN}{\overline\NN}
    \newcommand{\OO}{{\cal O}}
    \newcommand{\PP}{{\cal P}}
    \newcommand{\QQ}{{\cal Q}}
    \newcommand{\RR}{{\cal R}}
  \renewcommand{\SS}{{\cal S}}
    \newcommand{\TT}{{\cal T}}
    \newcommand{\uu}{{\cal u}}
    \newcommand{\UU}{{\cal U}}
    \newcommand{\VV}{{\cal V}}
    \newcommand{\XX}{{\cal X}}
    \newcommand{\YY}{{\cal Y}}
    \newcommand{\SOSO}{{\cal SO}}
    \newcommand{\GLGL}{{\cal GL}}

    \newcommand{\Ee}{{\rm E}}

  \newcommand{\NNN}{{\mathbb{N}}}
  \newcommand{\III}{{\mathbb{I}}}
  \newcommand{\ZZZ}{{\mathbb{Z}}}
  \newcommand{\RRR}{{\mathbb{R}}}
  \newcommand{\SSS}{{\mathbb{S}}}
  \newcommand{\CCC}{{\mathbb{C}}}
  \newcommand{\DDD}{{\mathbb{D}}}
  \newcommand{\one}{{{\bf 1}}}
  \newcommand{\eee}{\text{e}}

  \newcommand{\NNNN}{{\overline{\cal N}}}

  \renewcommand{\[}{\Big[}
  \renewcommand{\]}{\Big]}
  \renewcommand{\(}{\Big(}
  \renewcommand{\)}{\Big)}
  \renewcommand{\|}{\,|\,}
  \renewcommand{\;}{\,;\,}
  \renewcommand{\=}{\!=\!}
    \newcommand{\<}{\left\langle}
  \renewcommand{\>}{\right\rangle}

  \newcommand{\na}{{\nabla\!}}
  \newcommand{\he}{{\nabla^2\!}}
  \newcommand{\Prob}{{\rm Prob}}
  \newcommand{\Dir}{{\rm Dir}}
  \newcommand{\Beta}{{\rm Beta}}
  \newcommand{\Var}{{\rm Var}}
  \newcommand{\Aut}{{\rm Aut}}
  \newcommand{\cor}{{\rm cor}}
  \newcommand{\corr}{{\rm corr}}
  \newcommand{\cov}{{\rm cov}}
  \newcommand{\sd}{{\rm sd}}
  \newcommand{\tr}{{\rm tr}}
  \newcommand{\Tr}{{\rm Tr}}
  \newcommand{\rank}{{\rm rank}}
  \newcommand{\diag}{{\rm diag}}
  \newcommand{\id}{{\rm id}}
  \newcommand{\Id}{{\rm\bf I}}
  \newcommand{\Gl}{{\rm Gl}}
  \renewcommand{\th}{\ensuremath{{}^\text{th}} }
  \newcommand{\lag}{\mathcal{L}}
  \newcommand{\inn}{\rfloor}
  \newcommand{\lie}{\pounds}
  \newcommand{\longto}{\longrightarrow}
  \newcommand{\speer}{\parbox{0.4ex}{\raisebox{0.8ex}{$\nearrow$}}}
  \renewcommand{\dag}{ {}^\dagger }
  \newcommand{\blbox}{\rule{1ex}{1ex}}
  \newcommand{\Ji}{J^\sharp}
  \newcommand{\h}{{}^\star}
  \newcommand{\w}{\wedge}
  \newcommand{\too}{\longrightarrow}
  \newcommand{\oot}{\longleftarrow}
  \newcommand{\To}{\Rightarrow}
  \newcommand{\oT}{\Leftarrow}
  \newcommand{\oTo}{\Leftrightarrow}
  \renewcommand{\iff}{~\Longleftrightarrow~}
  \newcommand{\Too}{\;\Longrightarrow\;}
  \newcommand{\oto}{\leftrightarrow}
  \newcommand{\ot}{\leftarrow}
  \newcommand{\ootoo}{\longleftrightarrow}
  \newcommand{\ow}{\stackrel{\circ}\wedge}
  \newcommand{\feed}{\nonumber \\}
  \newcommand{\comma}{~,\quad}
  \newcommand{\period}{~.\quad}
  \newcommand{\del}{\partial}
  \newcommand{\point}{$\bullet~~$}
  \newcommand{\doubletilde}{ ~ \raisebox{0.3ex}{$\widetilde {}$} \raisebox{0.6ex}{$\widetilde {}$} \!\! }
  \newcommand{\topcirc}{\parbox{0ex}{~\raisebox{2.5ex}{${}^\circ$}}}
  \newcommand{\topdot} {\parbox{0ex}{~\raisebox{2.5ex}{$\cdot$}}}
  \newcommand{\topddot} {\parbox{0ex}{~\raisebox{1.3ex}{$\ddot{~}$}}}
  \newcommand{\sym}{\topcirc}
  \newcommand{\tsum}{\textstyle\sum}
  \newcommand{\st}{\quad\text{s.t.}\quad}

  \newcommand{\half}{\ensuremath{\frac{1}{2}}}
  \newcommand{\third}{\ensuremath{\frac{1}{3}}}
  \newcommand{\fourth}{\ensuremath{\frac{1}{4}}}

  \newcommand{\ubar}{\underline}
  \renewcommand{\vec}{\boldsymbol}
  \renewcommand{\*}{\text{\footnotesize\raisebox{-.4ex}{*}{}}}

  \newcommand{\gto}{{\raisebox{.5ex}{${}_\rightarrow$}}}
  \newcommand{\gfrom}{{\raisebox{.5ex}{${}_\leftarrow$}}}
  \newcommand{\gnto}{{\raisebox{.5ex}{${}_\nrightarrow$}}}
  \newcommand{\gnfrom}{{\raisebox{.5ex}{${}_\nleftarrow$}}}

  \DeclareMathOperator*{\argmax}{argmax}
  \DeclareMathOperator*{\argmin}{argmin}
  \DeclareMathOperator{\sign}{sign}
  \DeclareMathOperator{\acos}{acos}
  \newcommand{\ee}[1]{\ensuremath{\cdot10^{#1}}}
  \newcommand{\sub}[1]{\ensuremath{_{\text{#1}}}}
  \newcommand{\up}[1]{\ensuremath{^{\text{#1}}}}
  \newcommand{\kld}[3][{}]{D_{#1}\big(#2\,\big|\!\big|\,#3\big)}
  \newcommand{\sprod}[2]{\big<#1\,,\,#2\big>}
  \newcommand{\End}{\text{End}}
  \newcommand{\txt}[1]{\quad\text{#1}\quad}
  \newcommand{\Over}[2]{\genfrac{}{}{0pt}{0}{#1}{#2}}
  \newcommand{\arr}[2]{\hspace*{-.5ex}\begin{array}{#1}#2\end{array}\hspace*{-.5ex}}
  \newcommand{\mat}[3][.9]{
    \renewcommand{\arraystretch}{#1}{\scriptscriptstyle{\left(
      \hspace*{-1ex}\begin{array}{#2}#3\end{array}\hspace*{-1ex}
    \right)}}\renewcommand{\arraystretch}{1.2}
  }
  \newcommand{\Mat}[3][.9]{
    \renewcommand{\arraystretch}{#1}{\scriptscriptstyle{\left[
      \hspace*{-1ex}\begin{array}{#2}#3\end{array}\hspace*{-1ex}
    \right]}}\renewcommand{\arraystretch}{1.2}
  }
  \newcommand{\case}[2][ll]{\left\{\arr{#1}{#2}\right.}
  \newcommand{\seq}[1]{\textsf{\<#1\>}}
  \newcommand{\seqq}[1]{\textsf{#1}}
  \newcommand{\floor}[1]{\lfloor#1\rfloor}
  \newcommand{\Exp}[2]{\text{E}_{#1}\{#2\}}
  \newcommand{\ex}{\setminus}

  \providecommand{\href}[2]{{\color{blue}USE PDFLATEX!}}
  \providecommand{\url}[2]{\href{#1}{{\color{blue}#2}}}
  \newcommand{\anchor}[2]{\begin{picture}(0,0)\put(#1){#2}\end{picture}}
  \newcommand{\pagebox}{\begin{picture}(0,0)\put(-3,-23){
    \textcolor[rgb]{.5,1,.5}{\framebox[\textwidth]{\rule[-\textheight]{0pt}{0pt}}}}
    \end{picture}}

  \newcommand{\hide}[1]{
    \begin{list}{}{\leftmargin0ex \rightmargin0ex \topsep0ex \parsep0ex}
       \helvetica{5}{1}{m}{n}
       \renewcommand{\section}{\par SECTION: }
       \renewcommand{\subsection}{\par SUBSECTION: }
       \item[$~~\blacktriangleright$]
       #1%$\blacktriangleleft~~$
       \message{^^JHIDE--Warning!^^J}
    \end{list}
  }
  \newcommand{\Hide}{\renewcommand{\hide}[1]{\message{^^JHIDE--Warning (hidden)!^^J}}}
  \newcommand{\HIDE}{\renewcommand{\hide}[1]{}}
  \newcommand{\fullhide}[1]{}
  \newcommand{\todo}[1]{{\tt[TODO: #1]}\message{^^JTODO--Warning: #1^^J}}
  \newcommand{\Todo}{\renewcommand{\todo}[1]{\message{^^JTODO--Warning (hidden)!^^J}}}
  \newcommand{\myauthor}[1]{\author{#1}\newcommand{\theauthor}{#1}}%\@author}
  \newcommand{\mytitle}[1]{\title{#1}\newcommand{\thetitle}{#1}}%\@title}
  \newcommand{\header}{\begin{document}\mytitle\cleardefs}
  \newcommand{\contents}{{\tableofcontents}\renewcommand{\contents}{}}
  \newcommand{\footer}{\small\bibliography{marc,bibs}\end{document}}
  \newcommand{\widepaper}{\usepackage{geometry}\geometry{a4paper,hdivide={25mm,*,25mm},vdivide={25mm,*,25mm}}}
  \newcommand{\moviex}[2]{\movie[externalviewer]{#1}{#2}} %\pdflatex\usepackage{multimedia}
  \newcommand{\rbox}[1]{\fboxrule2mm\fcolorbox[rgb]{1,.85,.85}{1,.85,.85}{#1}}
  \newcommand{\mpage}[2]{{\begin{minipage}{#1\columnwidth}#2\end{minipage}}}
  \newcommand{\redbox}[2]{\fboxrule1mm\fcolorbox[rgb]{1,.7,.7}{1,.7,.7}{\begin{minipage}{#1\columnwidth}\center#2\end{minipage}}}
  \newcommand{\onecol}[2]{
    \begin{minipage}[c]{#1\columnwidth}#2\end{minipage}}
  \newcommand{\twocol}[5][0]{
    \begin{minipage}[c]{#2\columnwidth}#4\end{minipage}\hspace*{#1\columnwidth}%
    \begin{minipage}[c]{#3\columnwidth}#5\end{minipage}}
  \newcommand{\threecol}[7][0]{%
    \begin{minipage}[c]{#2\columnwidth}#5\end{minipage}\hspace*{#1\columnwidth}%
    \begin{minipage}[c]{#3\columnwidth}#6\end{minipage}\hspace*{#1\columnwidth}%
    \begin{minipage}[c]{#4\columnwidth}#7\end{minipage}}
  \newcommand{\threecoltext}[7][c]{
    \begin{minipage}[#1]{#2\textwidth}#5\end{minipage}%
    \begin{minipage}[#1]{#3\textwidth}#6\end{minipage}%
    \begin{minipage}[#1]{#4\textwidth}#7\end{minipage}}
  \newcommand{\threecoltop}[7][0]{%
   \begin{minipage}[t]{#2\columnwidth}#5\end{minipage}\hspace*{#1\columnwidth}%
   \begin{minipage}[t]{#3\columnwidth}#6\end{minipage}\hspace*{#1\columnwidth}%
   \begin{minipage}[t]{#4\columnwidth}#7\end{minipage}}
  \newcommand{\fourcol}[9][0]{%
   \begin{minipage}[c]{#2\columnwidth}#6\end{minipage}\hspace*{#1\columnwidth}%
   \begin{minipage}[c]{#3\columnwidth}#7\end{minipage}\hspace*{#1\columnwidth}%
   \begin{minipage}[c]{#4\columnwidth}#8\end{minipage}\hspace*{#1\columnwidth}%
   \begin{minipage}[c]{#5\columnwidth}#9\end{minipage}}
  \newcommand{\helvetica}[4]{\setlength{\unitlength}{1pt}\fontsize{#1}{#1}\linespread{#2}\usefont{OT1}{phv}{#3}{#4}}
  \newcommand{\helve}[1]{\helvetica{#1}{1.5}{m}{n}}
  \newcommand{\german}{\usepackage[german]{babel}\usepackage[utf8]{inputenc}}

\newcommand{\norm}[1]{|\!|#1|\!|}
\newcommand{\expr}[1]{[\hspace{-.2ex}[#1]\hspace{-.2ex}]}

\newcommand{\Jwi}{J^\sharp_W}
\newcommand{\THi}{T^\sharp_H}
\newcommand{\Jci}{J^\natural_C}
\newcommand{\hJi}{{\bar J}^\sharp}
\renewcommand{\|}{\,|\,}
\renewcommand{\=}{\!=\!}
\newcommand{\myminus}{{\hspace*{-.0pt}\text{\rm -}\hspace*{-.5pt}}}
\newcommand{\myplus}{{\hspace*{-.0pt}\text{\rm +}\hspace*{-.5pt}}}
\newcommand{\1}{{\myminus1}}
\newcommand{\2}{{\myminus2}}
\newcommand{\3}{{\myminus3}}
\newcommand{\mT}{{\text{\rm -}\hspace*{-1pt}\top}}
\newcommand{\po}{{\myplus1}}
\newcommand{\pt}{{\myplus2}}
\renewcommand{\T}{{\!\top\!}}
\newcommand{\xT}{{\underline x}}
\newcommand{\uT}{{\underline u}}
\newcommand{\zT}{{\underline z}}
\newcommand{\Sum}{\textstyle\sum}
\newcommand{\Int}{\textstyle\int}
\newcommand{\Prod}{\textstyle\prod}

\newenvironment{centy}{
\vspace{15mm}
\large
\hspace*{5mm}
\begin{minipage}{8cm}\it\color{blue}
}{
\end{minipage}
}

\newcommand{\old}{{\text{old}}}
\newcommand{\new}{{\text{new}}}
\newcommand{\MAP}{{\text{MAP}}}
\newcommand{\ML}{{\text{ML}}}

\newcommand{\redArrow}{\quad\anchor{0,-1}{\includegraphics[scale=.5]{figs/redArrow}}}
\newcommand{\pub}[1]{{\color{green}\helvetica{8}{1.3}{m}{n}#1\\}}
\DeclareMathOperator{\opKL}{KL}
\newcommand{\KL}[2]{\opKL\big(#1\,\big|\!\big|\,#2\big)} %\left(#1 |\!| #2\right)}

\renewcommand{\show}[2][.8]{\centerline{\includegraphics[width=#1\columnwidth]{#2}}}
\newcommand{\showh}[2][.8]{\includegraphics[width=#1\columnwidth]{#2}}
\newcommand{\shows}[2][.8]{\centerline{\includegraphics[scale=#1]{#2}}}
\newcommand{\showhs}[2][.8]{\includegraphics[scale=#1]{#2}}
\newcommand{\mov}[2]{\movie[externalviewer]{{\color{blue}\small #1}}{movies/#2}}
\newcommand{\movex}[2]{\movie[externalviewer]{#1}{#2}} %\pdflatex\usepackage{multimedia}
\newcommand{\cen}[1]{\centerline{#1}}

\newcommand{\redoMacrosInProof}{
  \renewcommand{\d}{\delta}
  \renewcommand{\=}{\!=\!}
}

\makeatletter
\newenvironment{code}{\smallskip\newline%
  \begin{lrbox}{\@tempboxa}\begin{minipage}{.99\columnwidth}\helvetica{7}{1.1}{m}{n}
}{
  \end{minipage}\end{lrbox}%
  \colorbox[rgb]{.95,.95,.95}{\usebox{\@tempboxa}}\smallskip\newline
}\makeatother

\newcommand{\refeq}[1]{(\ref{#1})}

\usepackage{algorithm}
\usepackage{algpseudocode}
\algrenewcommand{\algorithmicrequire}{\textbf{Input:~~}}
\algrenewcommand{\algorithmicensure}{\textbf{Output:}}
\algrenewcommand{\algorithmiccomment}[1]{\qquad\hfill~\hspace*{-5ex}\textit{// #1}}
\algrenewcommand{\alglinenumber}[1]{\helvetica{6}{1.3}{m}{n}#1:}

\newenvironment{algo}[1][8]{
\quad\begin{minipage}{.8\columnwidth}\helvetica{#1}{1.3}{m}{n}
\medskip\hrule\medskip
\begin{algorithmic}[1]
}{
\end{algorithmic}
\medskip\hrule\medskip
\end{minipage}
}

}

\newcommand{\stdpackages}{
  \usepackage{amsmath}
  \usepackage{amssymb}
  \usepackage{amsfonts}
  \allowdisplaybreaks
  \usepackage{amsthm}
  \usepackage{eucal}
  \usepackage{graphicx}
  \usepackage{color}
  \usepackage{geometry}

  \usepackage{multicol} 
  \usepackage{fancyhdr}

  \newcommand{\draft}{\usepackage[light,first]{draftcopy}\draftcopyName{draft}{350}}
  \newcommand{\labels}{\usepackage{showlabels}}
  \newcommand{\maple}{\usepackage{maple2e}}
  \newcommand{\makeidx}{\usepackage{makeidx}\makeindex}
  \newcommand{\chicago}{\usepackage{chicago}\bibliographystyle{chicago}
    \renewcommand{\refname}{References\renewcommand{\refname}{}}}
  \newcommand{\natbib}{\usepackage[round]{natbib}\bibliographystyle{abbrvnat}}
  \newcommand{\showlines}{
    \usepackage[modulo]{lineno} %options: pagewise, modulo, mathlines
    \renewcommand{\BM}{\begin{linenomath}}
    \renewcommand{\EM}{\end{linenomath}}
    \linenumbers
    \modulolinenumbers[5]
  }\newcommand{\BM}{}\newcommand{\EM}{}
}

\newcommand{\pdflatex}{
  \definecolor{bluecol}{rgb}{0,0,.5}
  \definecolor{greencol}{rgb}{0,.4,0}
  \usepackage[
    colorlinks,
    urlcolor=bluecol,
    citecolor=black,
    linkcolor=bluecol,
    pdfborder={0 0 0},
    pdfpagemode=UseNone, %UseOutlines, %bookmarks are displayed by acrobat
    pdfauthor={Marc Toussaint}
  ]{hyperref}
  \DeclareGraphicsExtensions{.pdf,.png,.jpg,.eps}
  \renewcommand{\r}{\varrho}
  \renewcommand{\l}{\lambda}
  \renewcommand{\L}{\Lambda}
  \renewcommand{\s}{\sigma}
  \renewcommand{\b}{\beta}
  \renewcommand{\d}{\delta}
  \renewcommand{\k}{\kappa}
  \renewcommand{\t}{\tau}
  \renewcommand{\O}{\Omega}
  \renewcommand{\o}{\omega}
  \renewcommand{\SS}{{\cal S}}
  \renewcommand{\=}{\!=\!}
}
\newcommand{\stdtheorems}{
  \theoremstyle{plain}
  \newtheorem{theorem}{Theorem}
  \newtheorem{lemma}[theorem]{Lemma}
  \newtheorem{corollary}[theorem]{Corollary}
  \newtheorem{proposition}{Proposition}
  \newtheorem{conjecture}{Conjecture}
  \newtheorem{result}{Result}[section]
  \newtheorem{hypothesis}{Hypothesis}[section]
  \theoremstyle{definition}
  \newtheorem{definition}{Definition}
  \theoremstyle{remark}
  \newtheorem{remark}{Remark}[section]
  \newtheorem{example}{Example}[section]
  \newtheorem{algoTheo}{Algorithm}
  \newtheorem{testTheo}{Test}
}
\newcommand{\stdstyle}[1]{
  \stdpackages
  \stdtheorems
  \renewcommand{\labelenumi}{\textbf{(\roman{enumi})}}
  \renewcommand{\theenumi}{(\roman{enumi})} %for ref
  \newcommand{\itemdot}{\renewcommand{\labelitemi}{\bf $\cdot$}}
  \newcommand{\enumA}{\renewcommand{\labelenumi}{\textbf{\Alph{enumi}}}}
  \newcommand{\blockindent}{3ex}
  \renewcommand{\baselinestretch}{#1}
  \renewcommand{\arraystretch}{1.2}
  \renewcommand{\topfraction}{1}
  \renewcommand{\bottomfraction}{1}
  \renewcommand{\textfraction}{0}
  \columnsep 5ex
  \parindent 3ex
  \parskip 1ex

  \parindent 0pt
  \topsep 4pt plus 1pt minus 2pt
  \partopsep 1pt plus 0.5pt minus 0.5pt
  \itemsep 2pt plus 1pt minus 0.5pt
  \parsep 2pt plus 1pt minus 0.5pt
  \parskip .5pc %add _in_ {thebibliography} environment in *.bbl

  \setcounter{tocdepth}{3}
  \setcounter{secnumdepth}{3}

  \geometry{a4paper,hdivide={35mm,*,35mm},vdivide={35mm,*,35mm}}

  \renewcommand{\headrulewidth}{.0pt}\renewcommand{\footrulewidth}{.0pt}\cfoot{}
  \fancyhead[OL,EC]{\it\theauthor---\today}
  \fancyhead[ER]{\leftmark}
  \fancyhead[OR,EL]{\thepage}
  \fancyfoot[EL,OR]{}
  \setlength{\headsep}{10mm}
  \renewenvironment{abstract}
    {\vspace*{5ex}\begin{rblock}\hrule\vspace{1.5ex}{\bf Abstract.~}\small}
    {\vspace{2ex}\hrule\end{rblock}\vspace{5ex}}
  \newenvironment{keyword}
    {\par{\it Keywords:~}}
    {}
  \newcommand{\published}{}
  \def\makemytitle{%
    \thispagestyle{empty}
    \begin{list}{}{\leftmargin3ex \rightmargin3ex \topsep0ex \parsep0ex}\item[]
      \begin{center}
        {\fontsize{18}{25}\selectfont{\thetitle\\}}\vspace{5ex}

        {\fontsize{14}{16}\selectfont{\theauthor\\}}\vspace{1ex}

        {\footnotesize{\sl \addressFUB}\\ \emailBerlin}

        {\footnotesize \today}

        \vspace{1ex}
        {\small \published}
      \end{center}
    \end{list}
    \renewcommand{\maketitle}{\chapter{\thetitle}}
  }
}

\newcommand{\cleardefs}{
  \renewcommand{\article}[2]{}
  \renewcommand{\book}[2]{}
  \renewcommand{\draft}{}
  \renewcommand{\labels}{}
  \renewcommand{\maple}{}
  \renewcommand{\makeidx}{}
  \renewcommand{\chicago}{}
  \renewcommand{\pdflatex}{}
  \renewcommand{\header}{}
}

\newcommand{\article}[2]{
  \stdstyle{#2}
  \usepackage{palatino}
  
}

\newcommand{\lectureNote}{
  \documentclass[10pt,fleqn,twocolumn]{article}
  \usepackage{amsmath}
  \usepackage{amssymb}
  \usepackage{amsfonts}
  \usepackage{amsthm}
  \usepackage{eucal}
  \usepackage{graphicx}
  \usepackage{color}
  \usepackage{fancyhdr}
  \usepackage{geometry}
  \usepackage{palatino}

  \renewcommand{\baselinestretch}{1.1}
  \geometry{a4paper,headsep=7mm,hdivide={15mm,*,15mm},vdivide={20mm,*,15mm}}
    
  \allowdisplaybreaks

  \fancyhead[OL,ER]{\thetitle, \textit{Marc Toussaint}---\today}
  \fancyhead[C]{}
  \fancyhead[OR,EL]{\thepage}
  \fancyfoot{}
  \pagestyle{fancy}
  
}

\newcommand{\slideScript}{
}

\newcommand{\nips}{
  \documentclass{article}
  \usepackage{nips07submit_e,times}
  \stdpackages
  \pagestyle{plain}
}

\newcommand{\nipsben}{
  \documentclass{article}
  \usepackage{nips06}
  \stdpackages
  \pagestyle{plain}
}

\newcommand{\ijcnn}{
  \documentclass[10pt,twocolumn]{ijcnn}
  \stdpackages
  \bibliographystyle{abbrv} 
}

\newcommand{\springer}{
  \documentclass{springer_llncs}
  \renewcommand{\theenumi}{\alph{enumi}}
  \renewcommand{\labelenumi}{(\alph{enumi})}
  \renewcommand{\labelitemi}{$\bullet$}
  \bibliographystyle{abbrv}
  \stdpackages\stdtheorems
}

\newcommand{\elsevier}{
  \documentclass{elsart1p}
  \usepackage{natbib}
  \bibliographystyle{elsart-harv}
  \stdpackages
  \stdtheorems

}

\newcommand{\ieeejournal}{
  \documentclass[journal,twoside]{IEEEtran}
  \renewcommand{\theenumi}{\roman{enumi}}
  \renewcommand{\labelenumi}{(\roman{enumi})}
  \bibliographystyle{IEEEtran.bst}
  \usepackage{cite}
  \stdpackages
  \stdtheorems
  
}

\newcommand{\ieeeconf}{
  \documentclass[a4paper, 10pt, conference]{ieeeconf}
  \IEEEoverridecommandlockouts
  \overrideIEEEmargins
  \bibliographystyle{IEEEtran.bst}
  \stdpackages
  \stdtheorems
  \renewcommand{\theenumi}{\roman{enumi}}
  \renewcommand{\labelenumi}{(\roman{enumi})}
  
}

\newcommand{\foga}{
  \documentclass{article} 
  \stdpackages
  \usepackage{foga-02}
  \usepackage{chicago}
  \bibliographystyle{foga-chicago}
}

\newcommand{\book}[2]{
  \documentclass[#1pt,twoside,fleqn]{book}
  \newenvironment{abstract}{\begin{rblock}{\bf Abstract.~}\small}{\end{rblock}}
  \stdstyle{#2}

}

\newcommand{\letter}{

  \stdstyle{1.1}
  \usepackage{palatino}
  \renewcommand{\a}{\alpha}
  \renewcommand{\b}{\beta}
  \renewcommand{\d}{\delta}
    \newcommand{\D}{\Delta}
    \newcommand{\e}{\epsilon}
    \newcommand{\g}{\gamma}
    \newcommand{\G}{\Gamma}
  \renewcommand{\l}{\lambda}
  \renewcommand{\L}{\Lambda}
    \newcommand{\m}{\mu}
    \newcommand{\n}{\nu}
    \newcommand{\N}{\nabla}
  \renewcommand{\k}{\kappa}
  \renewcommand{\o}{\omega}
  \renewcommand{\O}{\Omega}
    \newcommand{\p}{\varphi}
  \renewcommand{\P}{\Phi}
  \renewcommand{\r}{\varrho}
    \newcommand{\s}{\sigma}
  \renewcommand{\S}{\Sigma}
  \renewcommand{\t}{\theta}
    \newcommand{\T}{\Theta}
    \newcommand{\x}{\xi}
    \newcommand{\X}{\Xi}
    \newcommand{\Y}{\Upsilon}
    \newcommand{\z}{\zeta}

  \renewcommand{\AA}{{\cal A}}
    \newcommand{\BB}{{\cal B}}
    \newcommand{\CC}{{\cal C}}
    \newcommand{\cc}{{\cal c}}
    \newcommand{\DD}{{\cal D}}
    \newcommand{\EE}{{\cal E}}
    \newcommand{\FF}{{\cal F}}
    \newcommand{\GG}{{\cal G}}
    \newcommand{\HH}{{\cal H}}
    \newcommand{\II}{{\cal I}}
    \newcommand{\KK}{{\cal K}}
    \newcommand{\LL}{{\cal L}}
    \newcommand{\MM}{{\cal M}}
    \newcommand{\NN}{{\cal N}}
    \newcommand{\oNN}{\overline\NN}
    \newcommand{\OO}{{\cal O}}
    \newcommand{\PP}{{\cal P}}
    \newcommand{\QQ}{{\cal Q}}
    \newcommand{\RR}{{\cal R}}
  \renewcommand{\SS}{{\cal S}}
    \newcommand{\TT}{{\cal T}}
    \newcommand{\uu}{{\cal u}}
    \newcommand{\UU}{{\cal U}}
    \newcommand{\VV}{{\cal V}}
    \newcommand{\XX}{{\cal X}}
    \newcommand{\YY}{{\cal Y}}
    \newcommand{\SOSO}{{\cal SO}}
    \newcommand{\GLGL}{{\cal GL}}

    \newcommand{\Ee}{{\rm E}}

  \newcommand{\NNN}{{\mathbb{N}}}
  \newcommand{\III}{{\mathbb{I}}}
  \newcommand{\ZZZ}{{\mathbb{Z}}}
  \newcommand{\RRR}{{\mathbb{R}}}
  \newcommand{\SSS}{{\mathbb{S}}}
  \newcommand{\CCC}{{\mathbb{C}}}
  \newcommand{\DDD}{{\mathbb{D}}}
  \newcommand{\one}{{{\bf 1}}}
  \newcommand{\eee}{\text{e}}

  \newcommand{\NNNN}{{\overline{\cal N}}}

  \renewcommand{\[}{\Big[}
  \renewcommand{\]}{\Big]}
  \renewcommand{\(}{\Big(}
  \renewcommand{\)}{\Big)}
  \renewcommand{\|}{\,|\,}
  \renewcommand{\;}{\,;\,}
  \renewcommand{\=}{\!=\!}
    \newcommand{\<}{\left\langle}
  \renewcommand{\>}{\right\rangle}

  \newcommand{\na}{{\nabla\!}}
  \newcommand{\he}{{\nabla^2\!}}
  \newcommand{\Prob}{{\rm Prob}}
  \newcommand{\Dir}{{\rm Dir}}
  \newcommand{\Beta}{{\rm Beta}}
  \newcommand{\Var}{{\rm Var}}
  \newcommand{\Aut}{{\rm Aut}}
  \newcommand{\cor}{{\rm cor}}
  \newcommand{\corr}{{\rm corr}}
  \newcommand{\cov}{{\rm cov}}
  \newcommand{\sd}{{\rm sd}}
  \newcommand{\tr}{{\rm tr}}
  \newcommand{\Tr}{{\rm Tr}}
  \newcommand{\rank}{{\rm rank}}
  \newcommand{\diag}{{\rm diag}}
  \newcommand{\id}{{\rm id}}
  \newcommand{\Id}{{\rm\bf I}}
  \newcommand{\Gl}{{\rm Gl}}
  \renewcommand{\th}{\ensuremath{{}^\text{th}} }
  \newcommand{\lag}{\mathcal{L}}
  \newcommand{\inn}{\rfloor}
  \newcommand{\lie}{\pounds}
  \newcommand{\longto}{\longrightarrow}
  \newcommand{\speer}{\parbox{0.4ex}{\raisebox{0.8ex}{$\nearrow$}}}
  \renewcommand{\dag}{ {}^\dagger }
  \newcommand{\blbox}{\rule{1ex}{1ex}}
  \newcommand{\Ji}{J^\sharp}
  \newcommand{\h}{{}^\star}
  \newcommand{\w}{\wedge}
  \newcommand{\too}{\longrightarrow}
  \newcommand{\oot}{\longleftarrow}
  \newcommand{\To}{\Rightarrow}
  \newcommand{\oT}{\Leftarrow}
  \newcommand{\oTo}{\Leftrightarrow}
  \renewcommand{\iff}{~\Longleftrightarrow~}
  \newcommand{\Too}{\;\Longrightarrow\;}
  \newcommand{\oto}{\leftrightarrow}
  \newcommand{\ot}{\leftarrow}
  \newcommand{\ootoo}{\longleftrightarrow}
  \newcommand{\ow}{\stackrel{\circ}\wedge}
  \newcommand{\feed}{\nonumber \\}
  \newcommand{\comma}{~,\quad}
  \newcommand{\period}{~.\quad}
  \newcommand{\del}{\partial}
  \newcommand{\point}{$\bullet~~$}
  \newcommand{\doubletilde}{ ~ \raisebox{0.3ex}{$\widetilde {}$} \raisebox{0.6ex}{$\widetilde {}$} \!\! }
  \newcommand{\topcirc}{\parbox{0ex}{~\raisebox{2.5ex}{${}^\circ$}}}
  \newcommand{\topdot} {\parbox{0ex}{~\raisebox{2.5ex}{$\cdot$}}}
  \newcommand{\topddot} {\parbox{0ex}{~\raisebox{1.3ex}{$\ddot{~}$}}}
  \newcommand{\sym}{\topcirc}
  \newcommand{\tsum}{\textstyle\sum}
  \newcommand{\st}{\quad\text{s.t.}\quad}

  \newcommand{\half}{\ensuremath{\frac{1}{2}}}
  \newcommand{\third}{\ensuremath{\frac{1}{3}}}
  \newcommand{\fourth}{\ensuremath{\frac{1}{4}}}

  \newcommand{\ubar}{\underline}
  \renewcommand{\vec}{\boldsymbol}
  \renewcommand{\*}{\text{\footnotesize\raisebox{-.4ex}{*}{}}}

  \newcommand{\gto}{{\raisebox{.5ex}{${}_\rightarrow$}}}
  \newcommand{\gfrom}{{\raisebox{.5ex}{${}_\leftarrow$}}}
  \newcommand{\gnto}{{\raisebox{.5ex}{${}_\nrightarrow$}}}
  \newcommand{\gnfrom}{{\raisebox{.5ex}{${}_\nleftarrow$}}}

  \DeclareMathOperator*{\argmax}{argmax}
  \DeclareMathOperator*{\argmin}{argmin}
  \DeclareMathOperator{\sign}{sign}
  \DeclareMathOperator{\acos}{acos}
  \newcommand{\ee}[1]{\ensuremath{\cdot10^{#1}}}
  \newcommand{\sub}[1]{\ensuremath{_{\text{#1}}}}
  \newcommand{\up}[1]{\ensuremath{^{\text{#1}}}}
  \newcommand{\kld}[3][{}]{D_{#1}\big(#2\,\big|\!\big|\,#3\big)}
  \newcommand{\sprod}[2]{\big<#1\,,\,#2\big>}
  \newcommand{\End}{\text{End}}
  \newcommand{\txt}[1]{\quad\text{#1}\quad}
  \newcommand{\Over}[2]{\genfrac{}{}{0pt}{0}{#1}{#2}}
  \newcommand{\arr}[2]{\hspace*{-.5ex}\begin{array}{#1}#2\end{array}\hspace*{-.5ex}}
  \newcommand{\mat}[3][.9]{
    \renewcommand{\arraystretch}{#1}{\scriptscriptstyle{\left(
      \hspace*{-1ex}\begin{array}{#2}#3\end{array}\hspace*{-1ex}
    \right)}}\renewcommand{\arraystretch}{1.2}
  }
  \newcommand{\Mat}[3][.9]{
    \renewcommand{\arraystretch}{#1}{\scriptscriptstyle{\left[
      \hspace*{-1ex}\begin{array}{#2}#3\end{array}\hspace*{-1ex}
    \right]}}\renewcommand{\arraystretch}{1.2}
  }
  \newcommand{\case}[2][ll]{\left\{\arr{#1}{#2}\right.}
  \newcommand{\seq}[1]{\textsf{\<#1\>}}
  \newcommand{\seqq}[1]{\textsf{#1}}
  \newcommand{\floor}[1]{\lfloor#1\rfloor}
  \newcommand{\Exp}[2]{\text{E}_{#1}\{#2\}}
  \newcommand{\ex}{\setminus}

  \providecommand{\href}[2]{{\color{blue}USE PDFLATEX!}}
  \providecommand{\url}[2]{\href{#1}{{\color{blue}#2}}}
  \newcommand{\anchor}[2]{\begin{picture}(0,0)\put(#1){#2}\end{picture}}
  \newcommand{\pagebox}{\begin{picture}(0,0)\put(-3,-23){
    \textcolor[rgb]{.5,1,.5}{\framebox[\textwidth]{\rule[-\textheight]{0pt}{0pt}}}}
    \end{picture}}

  \newcommand{\hide}[1]{
    \begin{list}{}{\leftmargin0ex \rightmargin0ex \topsep0ex \parsep0ex}
       \helvetica{5}{1}{m}{n}
       \renewcommand{\section}{\par SECTION: }
       \renewcommand{\subsection}{\par SUBSECTION: }
       \item[$~~\blacktriangleright$]
       #1%$\blacktriangleleft~~$
       \message{^^JHIDE--Warning!^^J}
    \end{list}
  }
  \newcommand{\Hide}{\renewcommand{\hide}[1]{\message{^^JHIDE--Warning (hidden)!^^J}}}
  \newcommand{\HIDE}{\renewcommand{\hide}[1]{}}
  \newcommand{\fullhide}[1]{}
  \newcommand{\todo}[1]{{\tt[TODO: #1]}\message{^^JTODO--Warning: #1^^J}}
  \newcommand{\Todo}{\renewcommand{\todo}[1]{\message{^^JTODO--Warning (hidden)!^^J}}}
  \newcommand{\myauthor}[1]{\author{#1}\newcommand{\theauthor}{#1}}%\@author}
  \newcommand{\mytitle}[1]{\title{#1}\newcommand{\thetitle}{#1}}%\@title}
  \newcommand{\header}{\begin{document}\mytitle\cleardefs}
  \newcommand{\contents}{{\tableofcontents}\renewcommand{\contents}{}}
  \newcommand{\footer}{\small\bibliography{marc,bibs}\end{document}}
  \newcommand{\widepaper}{\usepackage{geometry}\geometry{a4paper,hdivide={25mm,*,25mm},vdivide={25mm,*,25mm}}}
  \newcommand{\moviex}[2]{\movie[externalviewer]{#1}{#2}} %\pdflatex\usepackage{multimedia}
  \newcommand{\rbox}[1]{\fboxrule2mm\fcolorbox[rgb]{1,.85,.85}{1,.85,.85}{#1}}
  \newcommand{\mpage}[2]{{\begin{minipage}{#1\columnwidth}#2\end{minipage}}}
  \newcommand{\redbox}[2]{\fboxrule1mm\fcolorbox[rgb]{1,.7,.7}{1,.7,.7}{\begin{minipage}{#1\columnwidth}\center#2\end{minipage}}}
  \newcommand{\onecol}[2]{
    \begin{minipage}[c]{#1\columnwidth}#2\end{minipage}}
  \newcommand{\twocol}[5][0]{
    \begin{minipage}[c]{#2\columnwidth}#4\end{minipage}\hspace*{#1\columnwidth}%
    \begin{minipage}[c]{#3\columnwidth}#5\end{minipage}}
  \newcommand{\threecol}[7][0]{%
    \begin{minipage}[c]{#2\columnwidth}#5\end{minipage}\hspace*{#1\columnwidth}%
    \begin{minipage}[c]{#3\columnwidth}#6\end{minipage}\hspace*{#1\columnwidth}%
    \begin{minipage}[c]{#4\columnwidth}#7\end{minipage}}
  \newcommand{\threecoltext}[7][c]{
    \begin{minipage}[#1]{#2\textwidth}#5\end{minipage}%
    \begin{minipage}[#1]{#3\textwidth}#6\end{minipage}%
    \begin{minipage}[#1]{#4\textwidth}#7\end{minipage}}
  \newcommand{\threecoltop}[7][0]{%
   \begin{minipage}[t]{#2\columnwidth}#5\end{minipage}\hspace*{#1\columnwidth}%
   \begin{minipage}[t]{#3\columnwidth}#6\end{minipage}\hspace*{#1\columnwidth}%
   \begin{minipage}[t]{#4\columnwidth}#7\end{minipage}}
  \newcommand{\fourcol}[9][0]{%
   \begin{minipage}[c]{#2\columnwidth}#6\end{minipage}\hspace*{#1\columnwidth}%
   \begin{minipage}[c]{#3\columnwidth}#7\end{minipage}\hspace*{#1\columnwidth}%
   \begin{minipage}[c]{#4\columnwidth}#8\end{minipage}\hspace*{#1\columnwidth}%
   \begin{minipage}[c]{#5\columnwidth}#9\end{minipage}}
  \newcommand{\helvetica}[4]{\setlength{\unitlength}{1pt}\fontsize{#1}{#1}\linespread{#2}\usefont{OT1}{phv}{#3}{#4}}
  \newcommand{\helve}[1]{\helvetica{#1}{1.5}{m}{n}}
  \newcommand{\german}{\usepackage[german]{babel}\usepackage[utf8]{inputenc}}

\newcommand{\norm}[1]{|\!|#1|\!|}
\newcommand{\expr}[1]{[\hspace{-.2ex}[#1]\hspace{-.2ex}]}

\newcommand{\Jwi}{J^\sharp_W}
\newcommand{\THi}{T^\sharp_H}
\newcommand{\Jci}{J^\natural_C}
\newcommand{\hJi}{{\bar J}^\sharp}
\renewcommand{\|}{\,|\,}
\renewcommand{\=}{\!=\!}
\newcommand{\myminus}{{\hspace*{-.0pt}\text{\rm -}\hspace*{-.5pt}}}
\newcommand{\myplus}{{\hspace*{-.0pt}\text{\rm +}\hspace*{-.5pt}}}
\newcommand{\1}{{\myminus1}}
\newcommand{\2}{{\myminus2}}
\newcommand{\3}{{\myminus3}}
\newcommand{\mT}{{\text{\rm -}\hspace*{-1pt}\top}}
\newcommand{\po}{{\myplus1}}
\newcommand{\pt}{{\myplus2}}
\renewcommand{\T}{{\!\top\!}}
\newcommand{\xT}{{\underline x}}
\newcommand{\uT}{{\underline u}}
\newcommand{\zT}{{\underline z}}
\newcommand{\Sum}{\textstyle\sum}
\newcommand{\Int}{\textstyle\int}
\newcommand{\Prod}{\textstyle\prod}

\newenvironment{centy}{
\vspace{15mm}
\large
\hspace*{5mm}
\begin{minipage}{8cm}\it\color{blue}
}{
\end{minipage}
}

\newcommand{\old}{{\text{old}}}
\newcommand{\new}{{\text{new}}}
\newcommand{\MAP}{{\text{MAP}}}
\newcommand{\ML}{{\text{ML}}}

\newcommand{\redArrow}{\quad\anchor{0,-1}{\includegraphics[scale=.5]{figs/redArrow}}}
\newcommand{\pub}[1]{{\color{green}\helvetica{8}{1.3}{m}{n}#1\\}}
\DeclareMathOperator{\opKL}{KL}
\newcommand{\KL}[2]{\opKL\big(#1\,\big|\!\big|\,#2\big)} %\left(#1 |\!| #2\right)}

\renewcommand{\show}[2][.8]{\centerline{\includegraphics[width=#1\columnwidth]{#2}}}
\newcommand{\showh}[2][.8]{\includegraphics[width=#1\columnwidth]{#2}}
\newcommand{\shows}[2][.8]{\centerline{\includegraphics[scale=#1]{#2}}}
\newcommand{\showhs}[2][.8]{\includegraphics[scale=#1]{#2}}
\newcommand{\mov}[2]{\movie[externalviewer]{{\color{blue}\small #1}}{movies/#2}}
\newcommand{\movex}[2]{\movie[externalviewer]{#1}{#2}} %\pdflatex\usepackage{multimedia}
\newcommand{\cen}[1]{\centerline{#1}}

\newcommand{\redoMacrosInProof}{
  \renewcommand{\d}{\delta}
  \renewcommand{\=}{\!=\!}
}

\makeatletter
\newenvironment{code}{\smallskip\newline%
  \begin{lrbox}{\@tempboxa}\begin{minipage}{.99\columnwidth}\helvetica{7}{1.1}{m}{n}
}{
  \end{minipage}\end{lrbox}%
  \colorbox[rgb]{.95,.95,.95}{\usebox{\@tempboxa}}\smallskip\newline
}\makeatother

\newcommand{\refeq}[1]{(\ref{#1})}

\usepackage{algorithm}
\usepackage{algpseudocode}
\algrenewcommand{\algorithmicrequire}{\textbf{Input:~~}}
\algrenewcommand{\algorithmicensure}{\textbf{Output:}}
\algrenewcommand{\algorithmiccomment}[1]{\qquad\hfill~\hspace*{-5ex}\textit{// #1}}
\algrenewcommand{\alglinenumber}[1]{\helvetica{6}{1.3}{m}{n}#1:}

\newenvironment{algo}[1][8]{
\quad\begin{minipage}{.8\columnwidth}\helvetica{#1}{1.3}{m}{n}
\medskip\hrule\medskip
\begin{algorithmic}[1]
}{
\end{algorithmic}
\medskip\hrule\medskip
\end{minipage}
}

  \parskip2.5ex
  \pagestyle{plain}
  \renewcommand{\familydefault}{\sfdefault}
}

\newcommand{\letterhead}[3]{
  \thispagestyle{empty}
  \vspace*{10mm}
  \begin{minipage}[t]{8cm}
    #1
  \end{minipage}
  \hspace*{\fill}
  \begin{minipage}[t]{6cm}
    \mbox{}~\hfill Prof.\ Dr.\ Marc Toussaint\\
    \mbox{}~\hfill Freie Universit\"at Berlin\\
    \mbox{}~\hfill Arnimallee 7\\
    \mbox{}~\hfill 14195 Berlin, Germany\\
    \mbox{}~\hfill +49 30 838 52485\\
    \mbox{}~\hfill marc-toussaint@fu-berlin.de
  \end{minipage}

   \vspace*{5mm}\hfill #3, \today\\

   \vspace*{5mm}{\textbf{#2}}

   \vspace*{5mm}
}

\newcommand{\slides}{
  \newcommand{\thepage}{\arabic{mypage}}
  \documentclass[t,hyperref={bookmarks=true}]{beamer}
  \usetheme{default}
  \usefonttheme[onlymath]{serif}
  \setbeamertemplate{navigation symbols}{}
  \setbeamersize{text margin left=5mm}
  \setbeamersize{text margin right=5mm}
  \setbeamertemplate{itemize items}{{\color{black}$\bullet$}}
  \graphicspath{{pics/}{figs/}}

  \stdpackages
  \usepackage{multimedia}

  \definecolor{bluecol}{rgb}{0,0,.5}
  \definecolor{greencol}{rgb}{0,.6,0}
  \renewcommand{\arraystretch}{1.2}
  \columnsep 0mm

  \columnseprule 0pt
  \parindent 0ex
  \parskip 0ex
  \newcommand{\headerfont}{\helvetica{14}{1.5}{b}{n}}
  \newcommand{\slidefont} {\helvetica{10}{1.4}{m}{n}}
  \renewcommand{\small} {\helvetica{9}{1.4}{m}{n}}
  \renewcommand{\tiny} {\helvetica{8}{1.3}{m}{n}}

  \newcounter{mypage}
  \newcommand{\incpage}{\addtocounter{mypage}{1}\setcounter{page}{\arabic{mypage}}}
  \setcounter{mypage}{0}
  \resetcounteronoverlays{page}

  \pagestyle{fancy}
  \renewcommand{\headrulewidth}{0pt} %1pt}
  \renewcommand{\footrulewidth}{0pt} %.5pt}
  \cfoot{}
  \rhead{}
  \lhead{}
  \rfoot{~\anchor{-10,12}{\tiny\textsf{\arabic{mypage}/\pageref{lastpage}}}}
  \definecolor{grey}{rgb}{.8,.8,.8}
  \definecolor{head}{rgb}{.85,.9,.9}
  \definecolor{blue}{rgb}{.0,.0,.5}
  \definecolor{green}{rgb}{.0,.5,.0}
  \definecolor{red}{rgb}{.8,.0,.0}
  \newcommand{\inverted}{
    \definecolor{main}{rgb}{1,1,1}
    \color{main}
    \pagecolor[rgb]{.3,.3,.3}
  }
  
}

\newcommand{\titleslide}[4][Marc Toussaint]{
  \newcommand{\slideauthor}{#1}
  \newcommand{\slidevenue}{#3}
  \slidefont
  \incpage
  \begin{frame}
  \begin{center}
    \vspace*{15mm}

    {\headerfont #2\\}
        
    \vspace*{7mm}

    #1 \\

    \vspace*{5mm}

    {\small 
      Machine Learning \& Robotics Lab -- University of Stuttgart\\
      marc.toussaint@informatik.uni-stuttgart.de

      \vspace*{3mm}

      \emph{#3}
    }

    \vspace*{0mm}

  \end{center}
  \begin{itemize}\item[]~\\
    #4
  \end{itemize}
  \end{frame}
}

\newcommand{\titleslideempty}[3]{
  \slidefont
  \incpage
  \begin{frame}
  \begin{center}
    \vspace*{15mm}

    {\headerfont #1\\}
        
    \vspace*{5mm}

    {\small\emph{#2}} \\

  \end{center}
  \begin{itemize}\item[]~\\
    #3
  \end{itemize}
  \end{frame}
}

\newcommand{\emptyslide}[2]{
  \slidefont
  \incpage
  \begin{frame}
    #2
  \end{frame}
}

\newcommand{\oldslide}[2]{
  \slidefont
  \incpage\begin{frame}
  \footskip-5mm
  \setlength{\unitlength}{1mm}
  {\headerfont #1} \vspace*{-2ex}
  \begin{itemize}\item[]~\\
    #2
  \end{itemize}
  \end{frame}
}

\newcommand{\slide}[2]{
  \slidefont
  \incpage\begin{frame}
  \addcontentsline{toc}{section}{#1}
  \vfill
  {\headerfont #1} \vspace*{-2ex}
  \begin{itemize}\item[]~\\
    #2
  \end{itemize}
  \vfill
  \end{frame}
}

\newcommand{\slidetop}[2]{
  \slidefont
  \incpage\begin{frame}
  ~
  {\headerfont #1} \vspace*{-2ex}
  \begin{itemize}\item[]~\\
    #2
  \end{itemize}
  \end{frame}
}

\newcommand{\slideempty}[2]{
  \slidefont
  \incpage\begin{frame}
  \setlength{\unitlength}{1mm}
  \begin{picture}(0,0)(6,4)
  \put(0,0){{\color{head}\rule{130mm}{13mm}{}}}
  \end{picture}
  {\headerfont #1}\\[-2ex]
    #2
  \end{frame}
}

\newcommand{\sliden}[2]{
  \slidefont
  \incpage\begin{frame}
  \setlength{\unitlength}{1mm}
  \centerline{\headerfont #1}
  \vspace*{-2ex}
  \begin{enumerate}\item[]~\\
    #2
  \end{enumerate}
  \end{frame}
}

\newcommand{\slidetwo}[2]{
  \slidefont
  \incpage\begin{frame}
  \setlength{\unitlength}{1mm}
  \begin{picture}(0,0)(6,5)
  \put(0,0){{\color{head}\rule{130mm}{14mm}{}}}
  \end{picture}
  {\headerfont #1}\\[-2ex]
  \begin{multicols}{2}
  \begin{itemize}%\item[]~\\
    #2
  \end{itemize}
  \end{multicols}
  \end{frame}
}

\newcommand{\poster}{
  \documentclass[fleqn]{article}
  \stdpackages
  \renewcommand{\baselinestretch}{1}
  \renewcommand{\arraystretch}{1.8}

  \usepackage{geometry}
  \geometry{
    paperwidth=1189mm,
    paperheight=841mm, %841mm, %91.3cm, % 120cm
    headheight=0mm,
    headsep=0mm,
    footskip=1mm,
    hdivide={3cm,*,3cm},vdivide={1cm,*,1cm}}

  \setlength{\columnsep}{5cm}
  \columnseprule 3pt
  \renewcommand{\labelitemi}{\rule[.4ex]{.6ex}{.6ex}~}

  \pagestyle{empty}

  \definecolor{grey}{rgb}{.9,.9,.9}
  \newcommand{\inverted}{
    \definecolor{main}{rgb}{1,1,1}
    \color{main}
    \pagecolor[rgb]{.3,.3,.3}
  }

}

\author{Marc Toussaint}

\newcommand{\inilogo}[1][.25]{\includegraphics[scale=#1]{INI}}
\newcommand{\rublogo}[1][.25]{\includegraphics[scale=#1]{RUB}}
\newcommand{\edinlogo}[1][.25]{\includegraphics[scale=#1]{pics/eushield-fullcolour}}

\newcommand{\addressCologne}{
  Institute for Theoretical Physics\\
  University of Cologne\\
  50923 K\"oln, Germany\\
  {\tt mt@thp.uni-koeln.de}\\
  {\tt www.thp.uni-koeln.de/\~{}mt/}
}

\newcommand{\emailINI}{mt@neuroinformatik.ruhr-uni-bochum.de}
\newcommand{\phoneINI}{+49-234-32-27974}
\newcommand{\faxINI}{+49-234-32-14209}
\newcommand{\urlINI}{\texttt{www.neuroinformatik.rub.de/PEOPLE/mt/}}
\newcommand{\AddressINI}{
  Institut~f\"ur~Neuroinformatik,
  Ruhr-Universit\"at~Bochum, ND~04,
  44780~Bochum, Germany
}
\newcommand{\addressINI}{
  Institut~f\"ur~Neuroinformatik\\
  Ruhr-Universit\"at Bochum, ND~04\\
  44780~Bochum, Germany
}

\newcommand{\emailANC}{mtoussai@inf.ed.ac.uk}
\newcommand{\phoneANC}{+44 131 650 3089}
\newcommand{\faxANC}{+44 131 650 6899}
\newcommand{\urlANC}{homepages.inf.ed.ac.uk/mtoussai}
\newcommand{\AddressANC}{
  School of Informatics,\\
  Institute for Adaptive and Neural Computation,\\
  University of Edinburgh, 5 Forrest Hill,\\
  Edinburgh EH1 2QL, Scotland, UK
}
\newcommand{\addressANC}{
  School~of~Informatics,\\
  University~of~Edinburgh, 5~Forrest~Hill\\
  Edinburgh~EH1~2QL, Scotland,~UK
}

\newcommand{\AddressBerlin}{
  Machine Learning \& Robotics group\\
  TU Berlin\\
  Franklinstr. 28/29, FR 6-9\\
  10587 Berlin, Germany
}
\newcommand{\addressTUB}{
  Machine~Learning~\&~Robotics~group, TU~Berlin\\\small
  Franklinstr. 28/29,~FR~6-9, 10587~Berlin, Germany
}
\newcommand{\addressFUB}{
  Machine~Learning~\&~Robotics~lab, FU~Berlin\\\small
  Arnimallee 7, 14195~Berlin, Germany
}
\newcommand{\addressUSTT}{
  Machine~Learning~\&~Robotics~lab, U~Stuttgart\\\small
  Universit{\"a}tsstra{\ss}e 38, 70569~Stuttgart, Germany
}
\newcommand{\emailBerlin}{mtoussai@cs.tu-berlin.de}
\newcommand{\phoneBerlin}{+49 30 314 24470}

\newcommand{\phoneHonda}{+49-69-89011-717}
\newcommand{\AddressHonda}{
  Honda Research Institute Europe\\
  Carl-Legien-Strasse 30\\
  63073 Offenbach/Main
}
\newcommand{\addressHonda}{
  Honda~Research~Institute~Europe~GmbH,\\\small
  Carl-Legien-Strasse~30, 63073~Offenbach/Main
}

\newlength{\subsecwidth}

\newcommand{\subsec}[1]{
  \addtocontents{toc}{
    \protect\setlength{\subsecwidth}{\textwidth}\protect\addtolength{\subsecwidth}{-27ex}
      \protect\vspace*{-1.5ex}\protect\hspace*{20ex}
      \protect\begin{minipage}[t]{\subsecwidth}\protect\footnotesize\protect\textsf{#1}\protect\end{minipage}
      \protect\par
  }
  \begin{rblock}\it #1\end{rblock}\medskip\noindent
}
\newcommand{\tocsep}{
  \addtocontents{toc}{\protect\bigskip}
}
\newcommand{\Chapter}[1]{
\chapter*{#1}\thispagestyle{empty}
\addcontentsline{toc}{chapter}{\protect\numberline{}#1}
}
\newcommand{\Section}[1]{
  \section*{#1}
  \addcontentsline{toc}{section}{\protect\numberline{}#1}
}
\newcommand{\Subsection}[1]{
  \subsection*{#1}
  \addcontentsline{toc}{subsection}{\protect\numberline{}#1}
}
\newcommand{\content}[1]{
}
\newcommand{\sepline}[1][200]{
  \begin{center} \begin{picture}(#1,0)
    \line(1,0){#1}
  \end{picture}\end{center}
}
\newcommand{\sepstar}{
  \begin{center} {\vspace{0.5ex}\rule[1.2ex]{5ex}{.1pt}~*~\rule[1.2ex]{5ex}{.1pt}} \end{center}\vspace{-1.5ex}\noindent
}
\newcommand{\partsection}[1]{
  \vspace{5ex}
  \centerline{\sc\LARGE #1}
  \addtocontents{toc}{\contentsline{section}{{\sc #1}}{}}
}
\newcommand{\intro}[1]{\textbf{#1}\index{#1}}

\newcounter{parac}
\newcommand{\para}{\noindent\refstepcounter{parac}{\bf [{\roman{parac}}]}~~}
\newcommand{\Pref}[1]{[\emph{\ref{#1}}\,]}

\newenvironment{itemS}{
\par
\tiny
\begin{list}{--}{\leftmargin4ex \rightmargin0ex \labelsep1ex
  \labelwidth2ex \topsep0pt \parsep0ex \itemsep0pt}
}{
\end{list}
}

\newenvironment{items}{
\par
\small%fontsize{9}{9}\linespread{1.2}
\begin{list}{--}{\leftmargin4ex \rightmargin0ex \labelsep1ex \labelwidth2ex
\topsep0pt \parsep0ex \itemsep3pt}
}{
\end{list}
}

\newenvironment{block}[1][]{{\noindent\bf #1}
\begin{list}{}{\leftmargin\blockindent \topsep-\parskip}
\item[]
}{
\end{list}
}

\newenvironment{rblock}{
\begin{list}{}{\leftmargin\blockindent \rightmargin\blockindent \topsep-\parskip}\item[]}{\end{list}}

\newcounter{algoi}
\newenvironment{algoList}{
\begin{list}{{(\thealgoi)}}
{\usecounter{algoi} \leftmargin7ex \rightmargin3ex \labelsep1ex
  \labelwidth5ex \topsep-.5ex \parsep.5ex \itemsep0pt}
}{
\end{list}\vspace*{1ex}
}
\newenvironment{test}[1][]{
  \medskip\begin{testTheo}[#1]~\begin{algoList}
}{
  \end{algoList}\end{testTheo}
}

\newcounter{questi}
\newenvironment{question}{
\begin{list}{\textbf{\thealgoi.}}
{\usecounter{algoi} \leftmargin2ex \rightmargin0ex \labelsep1ex
  \labelwidth1ex \topsep0ex \parsep.5ex \itemsep0pt}
\addtocounter{questi}{1}
\item[\textsf{Q\thequesti:}]
}{
\end{list}
}

\newenvironment{colpage}{
\addtolength{\columnwidth}{-3ex}
\begin{minipage}{\columnwidth}
\vspace{.5ex}
}{
\vspace{.5ex}
\end{minipage}
}

\newenvironment{enum}{
\begin{list}{}{\leftmargin3ex \topsep0ex \itemsep0ex}
\item[\labelenumi]
}{
\end{list}
}

\newenvironment{cramp}{
\begin{quote} \begin{picture}(0,0)
        \put(-5,0){\line(1,0){20}}
        \put(-5,0){\line(0,-1){20}}
\end{picture}
}{
\begin{picture}(0,0)
        \put(-5,5){\line(1,0){20}}
        \put(-5,5){\line(0,1){20}}
\end{picture} \end{quote}
}

\newcommand{\localcite}[1]{{
\begin{bibunit}[chicago]
\renewcommand{\refname}{\vspace{-\parskip}} \let\chapter\phantom \let\section\phantom
\nocite{#1}
\putbib[bibs]
\end{bibunit}%use \setcounter{enumiv}{xx} in thebibliography environment
}}

\newcommand{\boxpage}[2][\textwidth]{
  \setcounter{equation}{0}
  \renewcommand{\theequation}{A.\arabic{equation}}
  \fboxsep1ex
  \fbox{
  \begin{minipage}{#1}
    #2
  \end{minipage}
  }
}

\newcommand{\tightmath}{
  \setlength{\jot}{0pt}
  \setlength{\abovedisplayskip}{.5ex}
  \setlength{\belowdisplayskip}{.5ex}
}

\newenvironment{myproof}{
  \small
  \noindent \textit{Proof.~}
  \tightmath
}{
  \hfill
  \begin{picture}(0,0)(0,0)
  \put(-3.5,12){\rule{5pt}{5pt}}
  \end{picture}
}

\newcommand{\lst}{
  \usepackage{listings}
  \lstset{ %
    language=C,                % choose the language of the code
    basicstyle=\normalfont\small,       % the size of the fonts that are used for the code
    frame=none,                   % adds a frame around the code
    tabsize=4,                      % sets default tabsize to 2 spaces
    captionpos=b,                   % sets the caption-position to bottom
    texcl=true,
    mathescape=true,
    escapechar=\#,
    columns=flexible,
    xleftmargin=6ex,
    numbers=left, numberstyle=\footnotesize, stepnumber=1, numbersep=3ex
  }
}

  \stdstyle{1.1}
  \usepackage{palatino}

\pdflatex
\natbib

\newcommand{\comment}[1]{{\color{greencol}\helvetica{8}{1.2}#1\par}}
\renewcommand{\cen}{{\text{cen}}}
\newcommand{\any}{{\text{any}}}

\title{A Novel Augmented Lagrangian Approach for Inequalities and
Convergent Any-Time Non-Central Updates}
\date{October 1, 2014}

\author{Marc Toussaint\thanks{University of Stuttgart, Germany.
(\texttt{marc.toussaint@informatik.uni-stuttgart.de})}}

\begin{document}
\maketitle

\begin{abstract}
Motivated by robotic trajectory optimization problems we consider the
Augmented Lagrangian approach to constrained optimization.
%%  In
%% constrast 
%% The Augmented Lagrangian approach to constrained optimization has
%% recently received growing attention as an alternative to log-barrier
%% and SQP methods. We were motivated to revisit the approach for
%% robotics problems, where standard log-barrier methods compromise the
%% well-conditioning of the local Hessian of high-dimensional non-linear
%% trajectory optimization problems, which is avoided by the Augmented
%% Lagrangian approach.
We first propose an alternative augmentation of the Lagrangian to
handle the inequality case (not based on slack variables) and a
corresponding ``central'' update of the dual parameters. We proove
certain properties of this update: roughly, in the case of LPs and
when the ``constraint activity'' does not change between iterations,
the KKT conditions hold after just one iteration. This gives essential
insight on when the method is efficient in practise. We then present
our main contribution, which are consistent any-time (non-central)
updates of the dual parameters (i.e., updating the dual parameters
when we are not currently at an extremum of the Lagrangian). Similar
to the primal-dual Newton method, this leads to an algorithm that
parallely updates the primal and dual solutions, not distinguishing
between an outer loop to adapt the dual parameters and an inner loop
to minimize the Lagrangian. We again proof certain properties of this
anytime update: roughly, in the case of LPs and when constraint
activities would not change, the dual solution converges after one
iteration. Again, this gives essential insight in the caveats of the
method: if constraint activities change the method may destablize. We
propose simple smoothing, step-size adaptation and regularization
mechanisms to counteract this effect and guarantee monotone
convergence.
%%  We then extend the update to become 2nd order, leading to
%% the same properties also for (local) QPs. 
Finally, we evaluate the
proposed method on random LPs as well as on standard robot trajectory
optimization problems, confirming our motivation and intuition that
our approach performs well if the problem structure implies moderate
stability of constraint activity.
%% Based on this we define a concrete optimization algorithm \emph{Carla}
%% for non-linear constrained optimization
%% (with and without using 2nd order information in the dual updates) and
%% compare its performance on high-dimensional problems with that of
%% standard, non-concurrent Augmented Lagrangian, the interior point
%% method, and SQP.
\end{abstract}

\section{Introduction}

To motivate this work we first mention some empirical findings. We
tested standard interior point and Augmented Lagrangian methods on
random LPs and QPs as well as on non-linear constrained robot
trajectory optimization problems. For random LPs and QPs, we found
Augmented Lagrangian methods less efficient as plain
log-barrier. However, for our trajectory optimization problems
Augmented Lagrangian methods performed extremely well, only by a small
factor slower than unconstrained non-linear trajectory
optimization---in constrast to less efficient log-barrier methods. We
believe a reason for this is that in the trajectory optimization case
the constraints play a ``simpler'' role than in random LPs: the
problem is dominated by the non-linear cost function $f(x)$, the
number of constraints is smaller than the primal problem
dimensionality, and empirically we find that constraint activity is
rather stable, i.e., does not vary much over optimization interations.

These views motivate us to investigate in Augmented Lagrangian
methods, extending them to deal efficiently also with inequality
constraints particularly in cases where the constraint activity is
rather stable. We will propose an alternative augmentation to deal
with inequality, analyze it and generalize it towards an any-time
primal dual update. This analysis gives insight into why this
Augmented Lagrangian might be particularly appropriate when constraint
activity is rather stable. We will detail the contributions after
introducing related work.

%% constraints are heuristically translated to non-linear
%% but well-conditioned potentials, we can utilize highly efficient
%% iterative Newton methods to converge to local optima after only $\sim
%% 20$ iterations, leading to optimization times below one second and
%% online applicability in robotics. However, the heuristic

%% In these domains we optimize a trajectory $x_{0:T}$
%% where each $x_t \in\RRR^n$ with $n\approx 30$ and about $T\approx 200$
%% time slices. The cost function $f(x_{0:T})$ non-linear but usually
%% fairly smooth (e.g.\ involving squared potentials in trigonometric
%% functions of $x$) while the inequality constraint functions
%% $g(x_{0:T})$ can 

%% In the context of our concrete experience with high-dimensional
%% robotic trajectory optimization problems we got interested in the
%% Augmented Lagrangian approach because, unlike log-barrier and naive
%% squared penalty methods, it does not corrupt the well-conditioning of
%% the Hessian of the unconstrained problem. We found this to be a
%% crucial advantage, as taking the limit $\mu\to 0$ in log-barrier
%% methods, for instance, may degrade the efficiency of iterative Newton
%% methods for centering.

\section{Related work}

In Section 17.4, \cite{nocedal1999numerical} propose an ``unconstrained
formulation'' of the Augmented Lagrangian in the case of
inequalities. The specific dual parameter update (their Eq.\ (17.63))
is the same as the ``central'' update we consider below. However,
their specific augmentation (17.64) is different to the Augmented
Lagrangian we will propose---only includes a squared penalty $g^2$ if
$2\mu g+\l \le0$ (translated to our notation). We will explicitly
address the difference when discussing the implications of our
choice. \cite{nocedal1999numerical} state that their proposition has
not been practically evaluated and we are not aware of evaluations of
their approach. Futher, they do not extend towards any-time primal-dual
updates.% and 2nd order updates.

LANCELOT is the most popular software using the Augmented Lagrangian
for globally convergent non-linear
optimization \citep{conn1991globally,conn2010lancelot}. Inequalities
are handled with slack variables $\xi$, which implies that the
dimensionality of the optimization problem is increased and the state
space will be subject to bound constraints ($\xi\ge
0)$ \citep{nocedal1999numerical}, prohibiting straight-forward Newton
methods. Both of these aspects makes the approach less attractive in
the high-dimensional trajectory optimization domain. Further, we are
not aware of any-time updates used within such approaches.

Another approach is to consider shifted barriers (e.g., log-barriers)
as augmentation in the inequality
case \citep{conn1997globally,noll2004partially,noll2007local}. We find
these approaches very interesting and at first sight very different to our
$[\l>0 \vee g>0] g^2$ augmentation we will discuss below.
Again, we are not aware of any-time updates having been proposed for
such types of augmentations. We believe our approach to any-time
updates could be generalized also to the case of shifted barrier
augmentations.

We would also like to point to a very interesting historical
discussion of interior point methods by \cite{forsgren2002interior},
where the authors nicely clarify the original motivation for Augmented
Lagrangian methods: Log-barrier (and squared penalty) methods lead to
an ill-conditioning of the Hessian in the limit of $\m\to 0$ (strict
barriers). This was considered a problem and motivation for the
Augmented Lagrangian, which happens to beautifully not modify the
conditioning of the Hessian at all. However, in the late 80ies it was
thoroughly understood that the log-barrier's ill-conditioning of the
Hessian is, surprisingly, not a problem (confirming the practical
success), which lead to the rise of interior point methods and
efficient primal-dual formulations, diminishing the interest in the
Augmented Lagrangian. As mentioned in the introduction we feel that it
very much depends on the concrete structure of the problem whether
interior point or Augmented Lagrangian methods might be more
efficient.

 %vanderbei1999interior
%% \cite{goldfarb1999modified}

%% \item Ben taskar extra gradients [Nathan]

%% [[move]] Example: $\l$ is very large; we get below the constraint $g<0$; still
%% the penalty pulls it to the constraint; just as for equality that
%% allows to adjust $\l$ correctly. Without the pulling to the
%% constraint---although $g<0$---we could not calibrate $\l$ better. And,
%% as $\l$ is very large, the solution closer to the constraint is
%% perhaps also the more correct one (if $\l$ turns out to be active
%% eventually).

Our contributions over this previous work are:

\noindent (1) We propose an alternative augmentation for the inequality
case. We analyze the properties of a centered update (also proposed
in \citep[Eq.~(17.63)]{nocedal1999numerical}) with this augmentation,
giving sufficient conditions for when the update yields the dual
solution. This result gives essential insights on when the approach is
promising in practise.

\noindent (2) Based on these results we reason about which any-time
(non-centered) dual update (i.e., an update of dual parameters while
not being at an extremum of the Lagrangian) would also yield the
correct dual solution (under similar sufficient conditions). We
propose such an any-time update and provide these sufficient
conditions, generalizing the result of the first part.

\noindent (3) Finally we consider a straight-forward extension to
account for the local Hessian of $f(x)$, leading to a 2nd order
any-time update.

\section{Alternative augmentation and centered update}

Let $x\in\RRR^n$, $f:~ \RRR^n \to \RRR$, $g:~ \RRR^n \to \RRR^m$,
$h:~ \RRR^n \to \RRR^l$. We consider
\begin{align}
\min_x~ f(x) \st g(x)\le 0,~ h(x) = 0 ~.
\end{align}
We denote the dual variables as $\l\in\RRR^m, \k\in\RRR^l$. The KKT
conditions are
\begin{align}
\na f(x) + \l^\T \na g(x) + \k^\T \na h(x) &= 0 && \text{(stationarity)}\\
g(x) \le 0 ~\wedge~ h(x)&=0 && \text{(primal feasibility)}\\
\l &\ge 0  && \text{(dual feasibility)}\\
\l \cdot g_i(x) &= 0 && \text{(complementary)}
\end{align}
Let use introduce some notation we use throughout. By
primal-dual \emph{state} we refer to an arbitrary tuple
$(x,\l,\k)$. In any state $(x,\l,\k)$ we call the $i$th
constraint \emph{active} iff $\l_i>0 \vee g_i(x)>0$. For two vectors
$v$ and $w$, $(v;w)=(v^\T,w^\T)^\T$ denotes their ``stacking''
(analogously for matrices).

We consider the following Augmented Lagrangian,
which includes for any active constraint a squared penalty $g_i(x)^2$
pulling $g_i$ to zero.
\begin{definition} We define our Augmented Lagrangian as
\begin{align}\label{eqFaula}
%&\min_x~ L(x,\l,\k,f_0) ~, \\
L(x,\l,\k)
&= f(x)
 + \mu \sum_{i=1}^m [\l_i>0 \vee g_i(x)>0]~ g_i(x)^2
 + \l^\T g(x)  \feed
&~\phantom{= f(x)}
 + \nu \sum_{i=1}^l h_i(x)^2 
 + \k^\T h(x) \\
&= f(x)
 + [\mu I_\l(x) g(x) + \l]^\T g(x)
 + [\nu h(x) + \k]^\T h(x) ~,
\end{align}
where $I_\l(x) := \diag([g_i(x)\ge 0 \vee \l_i>0])$.
\end{definition}
Its gradient is
\begin{align}
\na L(x,\l,\k)
&= \na f(x)
 + [2 \mu I_\l g(x) + \l]^\T \na g(x)
 + [2 \nu h(x) + \k]^\T \na h(x) ~,
\end{align}
\begin{definition}
For any state $(x,\l,\k)$ we define the \emph{centered} update
$\UU^\cen$ as
\begin{align}
\UU^\text{cen}(x,\l,\k)
&= (\l',\k')\quad \text{with} \\
\l'
&= \max\{0, \l + 2 \mu g(x)\}  \label{eqLambdaUpdate}\\
\k'
&= \k + 2 \nu h(x) ~,
\end{align}
where the $\max{}$ operator is interpreted element-wise.
\end{definition}
This update is also introduced
in \citep[Eq.~(17.63)]{nocedal1999numerical}).  The centered update is
meant to be applied at a minimum $x'= \argmin_x~ L(x,\l,\k)$. The
standard nested loop approach uses an inner loop to converge to
$x'= \argmin_x~ L(x,\l,\k)$ for given dual parameters and an outer
loop to update $(\l',\k') = \UU^\cen(x',\l,\k)$.  An intuition behind
the update, related to the following theorem, is the following:
Assuming initially $\l=\k=0$, $x'=\argmin_x L(x,0,0)$ will
violate constraints. The squared penalties counteract these violations
by generating the gradients $2 \mu g(x')^\T I_\l(x') \na g_i(x) +
2 \nu h(x')^\T \na h(x')$ at $x'$. The centered update will generate exactly
these gradients in the next iteration. In other terms, the dual
parameters are chosen such that $\l^\T g(x) + \k^\T h(x)$ will
generated the gradients that have previously been generated by the
squared penalties. This is made more rigorous in the following result.

\begin{theorem}
For any $(\l,\k)$, let
\begin{align}
x'=\argmin_x L(x,\l,\k) \comma
(\l',\k')=\UU^\cen(x',\l,\k) \comma
x^* = \argmin_x L(x,\l',\k') ~.
\end{align}
For any Linear Program ($f$, $g$ and $h$ linear), if all active
constraints are linearly independent at $x^*$ (non-zero rows of
$I_{\l'}(x^*)~ \na g(x^*)$ are linearly independent), then
\begin{align}
[\forall_i: \l_i>0 \To \l_i'>0] ~\To~ \text{KKT hold at $x^*$}
\end{align}
\end{theorem}
\begin{proof}
Note that (element-wise)
\begin{align}\label{condition}
[\l>0 \To \l'>0 ]
&~\Longleftrightarrow~ [\l=0 \vee \l'>0] \\
&~~\Longrightarrow~ [\max\{0, \l + 2\mu g(x')\} = \l + 2 \mu I_\l(x') g(x')] ~.
\end{align}
This is obvious for $\l=0$. In the case $\l>0 \wedge \l'>0$ we have
$I_\l(x')=1$ and $\l+2\mu g(x') > 0$, from which the RHS follows.

We consider the gradient at $x'$,
\begin{align}\label{eq0}
0
&= \na L(x',\l,\k)
 = \na f
 + [2 \mu I_\l(x') g(x') + \l]^\T \na g
 + [2 \nu h(x') + \k]^\T \na h ~,
\end{align}
where $\na f, \na g, \na h$ are independent of $x'$, 
and compare it to the gradient at $x^*$,
\begin{align}
0
&= \na L(x^*, \l', \k') \\
&= \na f
 + [2 \mu I_{\l'}(x^*) g(x^*) + \max\{0, \l + 2 \mu g(x')\}]^\T \na g
 + [2 \nu h(x^*) + \k + 2 \nu h(x')]^\T \na h \label{eq3}\\
&= \na f
 + [2 \mu I_{\l'}(x^*) g(x^*) + \l + 2 \mu I_\l(x') g(x')]^\T \na g
 + [2 \nu h(x^*) + \k + 2 \nu h(x')]^\T \na h ~, \label{eq1} \\
&= 2 \mu g(x^*)^\T I_{\l'}(x^*) \na g
 + 2 \nu h(x^*)^\T \na h ~,
\end{align}
where in the 3rd line we use the implication of
$[\forall_i: \l_i>0 \To \l_i'>0 ]$, and the last line
inserted \refeq{eq0}. If all non-zero rows of $I_{\l'}(x^*) \na g$ are
linearly independent the gradient $\na L(x^*,\l',\k')$ is zero at
$x^*$ iff
\begin{align}\label{eq2}
I_{\l'}(x^*) g(x^*)=0,~ h(x^*)=0 ~.
\end{align}
Note that \refeq{eq2} implies primal feasibility, complementarity, as
well as $\UU(x^*,\l',\k') = (\l',k')$. With $\na L(x^*,\l',\k')=0$
and \refeq{eq2} the stationarity holds. Dual feasibility holds
by construction.
\end{proof}
The theorem states that, for a linear program, the updated dual
parameters $(\l',\k')$ are optimal under two conditions: 1)
$[\l>0 \To \l'>0]$, that is, none of the constraints becomes inactive
when it was previously active. And 2), all active constraints
(non-zero rows of $I_{\l'}(x^*) \na g$) are linearly independent. The
discussion of these two conditions is interesting and gives insight
into our choice of the augmentation itself.

Let us first discuss the case $[\l>0 \To \l'>0]$, where in some
iteration $\l>0$, then $x'=\argmin_x L(x,\l,\k)$ pushes far outside
the constraint ($2 \mu g(x')<-\l<0$) such that the subsequent update
chooses $\l'=0$. In this case, even in the locally linearized view,
the next optimization does not reach a KKT point. This is intuitive as
the inner loop optimization of $x'$ considered the constraint to be
strictly active and therefore included a penalty $g(x)^2$ \emph{even}
when $g(x)<0$. It pulled towards the constraint $g(x)=0$ even though
$g(x)<0$. This explains that the update failed to lead to a KKT point
directly: In the equations we see that \refeq{eq3} becomes unequal
to \refeq{eq1} when the $\max{}$ selects $\l'=0$ while $\l>0$, and
therefore the $\l'$ does not generate the necessary gradients to
achieve stationarity in the next centering. A trivial solution seems
to initialize $\l=0$ in the first iteration, which avoids
$[\l>0 \To \l'>0]$; however, here the second conditions gets into
play.

Whether non-zero rows of $I_{\l'}(x^*) \na g$ are linear independent
typically depends on $\sum_i [\l_i'>0 \vee g_i(x^*)>0] \le n$, that
is, how many constraints are ``active''. Note that $\l'$ has been
computed at $x'$ while $g(x^*)$ is evaluated at $x^*$. Therefore,
$I_{\l'}(x^*)$ includes constraints that have been active at
$x'$ \emph{or} at $x^*$, which can well be more than $n$.

If both conditions are fulfilled, the Theorem shows it is effective to
include the penalty $g(x)^2$ \emph{even} when $g(x)<0$, because it
leads to a `correct' retuning of the active dual parameter
$\l'>0$---under the given assumptions. It penalizes the inequality
just like an equality, assuming that $\l$ might remain active when it
was active before.

\section{Any-time, non-centered update}

\begin{definition}We define the any-time update as
\begin{align}
\UU^\any(x,\l,\k)
&= (\l', \k') \quad\text{with} \\
\mat{c}{\l'\\ \k'}
&= \argmin_{(\hat\l;\hat\k):\hat\l\ge 0} \bigg|\!\bigg|[\mat{c}{\hat\l\\ \hat\k}-\mat{c}{\l + 2 \mu I_\l(x) g(x) \\ \k +
2 \nu h(x)}]^\T \mat{c}{\na g(x)\\ \na h(x)} + \na L(x,\l,\k) \bigg|\!\bigg|^2
\label{eqAnyUp}
\end{align}
\end{definition}
This update is a bounded quadratic program, aiming to minimize the
difference between the gradients $[\l + 2 \mu I_\l(x') g(x')]^\T \na g
+ [\k + 2 \nu h(x')]^\T \na h - \na L(x',\l,\k)$ before the update,
and $\l'^\T \na g  + \k'^\T \na h$ after the update.

\begin{theorem}\label{thmMain}
Let $(x',\l,\k)$ be arbitrary and
\begin{align}
(\l',\k')=\UU^\any(x',\l,\k) \comma
x^* = \argmin_x L(x,\l',\k') ~.
\end{align}
For any Linear Program ($f$, $g$ and $h$ linear), if all active
constraints are linearly independent at $x^*$ (non-zero rows of
$I_{\l'}(x^*)~ \na g(x^*)$ are linearly independent), and if the
$\argmin$ in the update \refeq{eqAnyUp} reaches zero, then KKT hold at $x^*$.
\end{theorem}

\begin{proof}
We have
\begin{align}
\na L(x,\l',\k')
&=  \na f
 + [2 \mu I_{\l'}(x) g(x)+\l']^\T \na g
 + [2 \nu h(x)+\k']^\T \na h \\
&= \na f
 + [2 \mu I_{\l'}(x) g(x) + \l + 2 \mu I_\l(x') g(x')]^\T \na g \feed
&~\phantom{=\na f}
 + [2 \nu h(x) + \k + 2 \nu h(x')]^\T \na h ~ - ~ \na L(x',\l,\k)\\
&= 2 \mu g(x)^\T I_{\l'}(x) \na g
 + 2 \nu h(x)^\T \na h ~,
\end{align}
The rest of the proof is as previously.
\end{proof}

The above theorem makes a statement under the strong assumption
that we can minimize the $\argmin$ in \refeq{eqAnyUp} to zero. A
particular complication here is the bound constraint $\l'\ge 0$ of the
minimization, which in the centered update translated to the
$\l'\gets \min\{...,0\}$, which in turn was related to the
assumption $\l_i>0 \To \l_i'>0$ we made in Theorem 1.

To avoid the bounded optimization problem \refeq{eqAnyUp} we first
consider even stronger assumption which leads an analytical solution:
\begin{corollary}
Let $\forall_i:~ \l_i>0 \oTo \l'_i>0 \oTo g_i(x')\ge 0$, then the
minimum of \refeq{eqAnyUp} is given analytically as
\begin{align}
\mat{c}{\l'\\ \k'}
&= y - (A A^\T)^\1 A \na L(x,\l,\k) \comma
y = \mat{c}{\l + 2 \mu g(x) \\ \k + 2 \nu h(x)}_{>0} \comma
A = \mat{c}{\na g(x)\\ \na h(x)}_{>0} ~,
\end{align}
where the notation $(~)_{>0}$ refers to rows for which $\l_i>0$ only.
\end{corollary}
\begin{proof}

Under the strong assumption, all inactive constraints drop out of
the minimization \refeq{eqAnyUp} (as when and (as $\l'_i>0$ for the active
ones) \refeq{eqAnyUp} becomes and unconstrained minimization that can
be solve analytically. We have
\begin{align}
\mat{c}{\l'\\ \k'}
&= \argmin_{(\hat\l;\hat\k)} \bigg|\!\bigg| [\mat{c}{\hat\l\\ \hat\k}-\underbrace{\mat{c}{\l + 2 \mu g(x) \\ \k +
2 \nu h(x)}}_y]^\T \underbrace{\mat{c}{\na g(x)\\ \na h(x)}}_A + \na
L(x,\l,\k)  \bigg|\!\bigg|^2 \\
&= y - (A A^\T)^\1 A \na L(x,\l,\k) ~, \label{eq4}
\end{align}
which gives the minimum via a pseudo-inverse of the active constraint
matrix $A$.
\end{proof}

In our evaluations we employed an approximation to \refeq{eqAnyUp},
where we analytically solve the unconstrained problem \refeq{eq4} and
then impose the bound $\l'\ge 0$ by clipping values.

\begin{corollary}
$$\na L(x,\l,\k)=0 \quad\To\quad \UU^\any(x,\l,\k) = \UU^\cen(x,\l,\k)$$
\end{corollary}
That is, when $\na L(x,\l,\k)=0$ the any-time update coincides with the
centered update---as the argmin reaches zero when for $(\l', \k') =
 \UU^\cen(x,\l,\k)$.

\paragraph{Heuristic update.} Solving the bound constraint
problem \refeq{eqAnyUp} becomes yet another constrained optimization
problem. A heuristic is to update with \refeq{eq4} and then truncate
$\l_i' \gets \max\{0, \l_i'\}$. Again, for $\na L=0$ this coincides
with the centered update. For $\na L\not=0$ this is clearly a
suboptimal update. Empirical studies need to evaluate the benefit of
the any-time update.

%% \section{2nd order updates}

%% All of the above should also translate to having a 2nd order $f$ (but
%% linear constraints).

%% \begin{align}
%% L(x)
%%  &= f(x)
%%   + \mu \sum_{i=1}^m [g_i(x)\ge 0 \vee \l_i>0]~ g_i(x)^2
%%   + \nu \sum_{i=1}^l h_i(x)^2 + \l^\T g(x) + \k^\T h(x) \\
%% \na L(x)
%%  &= \na f
%%  + \sum_{i=1}^m [2\mu [..] g_i + \l_i] \na g_i
%%  + \sum_{i=1}^l [2\nu [..] h_i + \k_i] \na h_i \\
%% \he L(x)
%%  &\approx \he f
%%  + \sum_{i=1}^m 2\mu [..] \na g_i \na g_i^\T
%%  + \sum_{i=1}^l 2\nu [..] \na h_i \na h_i^\T
%% \end{align}

%% Newton steps use the inverse hessian -- what do we learn??

\section{Experiments}

\subsection{Algorithmic details}

In our experiments we use a basic Newton method for solving the
unconstrained problem $x'= \argmin_x~ L(x,\l,\k)$ up to a stopping
criterion. The method includes adaptive stepsize and
Levenberg-Marquardt damping, see Algorithm \ref{NA}. In the case of
the any-time update we increase the tolerance $\d$ by a factor 2 in
each Newton step, leading to an early stopping such that
$x' \approx \argmin_x~ L(x,\l,\k)$ only crudely approximates the
Lagrangian minimum. This is then alternated with the any-time
upate. To ensure that the monotonicity check (line 8) remains
sensible, the any-time update also needs to update the stored values
of $(L, \na L, \he L)$ (stored in line 9) consistently.

\begin{algorithm}[t]
\caption{\label{NA}Newton with adaptive step size and
Levenberg-Marquardt parameter}
\begin{algorithmic}[1]\small
\Require start point $x$, tolerance $\d$, functions $x \mapsto
  (f(x), \na f(x), \he f(x))$, parameters (defaults:
  $\a^0=\b^0=1, \a^+=2, \a^-=0.1, \b^+=\b^-=1, \r=0.01$)
\Ensure converged point $x$
\State initialize $\a=\a^0$, $\b=\b^0$
\State compute $f,\na f,\he f$ at $x$
\Repeat
\State compute $\D$ to solve $(\he f + \b \Id)~ \D = - \na f$
\Repeat \Comment{backtracking line search (for $\b^+=1$)}
\State $x' \gets x + \a\D$
\State compute $f',\na f', \he f'$ at $x'$ \Comment{computing $\na
f, \he f$ can be postponed}
\If{$f'\le f + \r\a \na f(x)^\T \D$} \Comment{step is accepted (Wolfe condition)}
\State $x \gets x',\quad (f,\na f,\he f) \gets (f',\na f', \he f')$
\State $\b \gets \b^-\b,\quad \a \gets \min\{\a^+\a,1\}$ \Comment{adapt $\a$
towards 1}
\Else \Comment{step is rejected}
\State\textbf{if} $\a\norm{\D}_\infty\ll\d$ or evaluations exceed \textbf{then}
abort with failure \Comment{gradient seems incorrect}
\State $\b \gets \b^+\b,\quad \a \gets \a^-\a$
\EndIf
\Until step accepted or $\b^+\not=1$ \Comment{change of $\b$
requires recomputing $\D$}
\Until $\b\le 1 \wedge \norm{\D}_\infty < \d$ or evaluations exceed
\end{algorithmic}
\end{algorithm}

\subsection{Random LPs}

\begin{figure}[t]
\label{figRandLPs}
\show[.6]{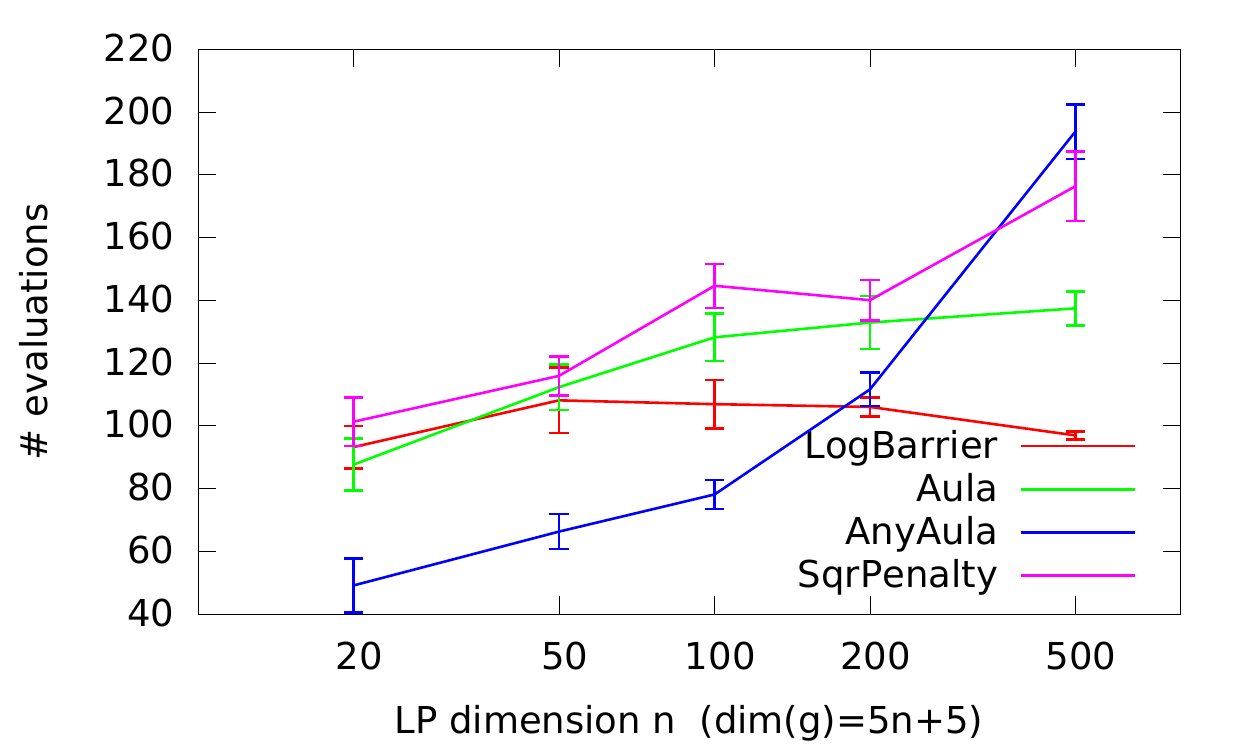}
\caption{Number of evaluations until convergence (tolerance $10^{-4}$)
  averaged over 10 random LPs for different dimensions
  $n=\dim(x)$. Errorbars indicate the deviation of the mean estimator.}
\end{figure}

We first compare the performance on random $n$-dimensional LPs of the form
\begin{align}
\min_x \sum_{i=1}^n x_i \st G~ \mat{c}{1\\x} \le 0
\end{align}
where the constraint-defining matrix $G\in\RRR^{m\times n\po}$ was
randomly generated as follows: First, each $G_{ij}\sim\NN(0,1)$;
second, if $G_{i1}>0:~ G_{i1}\gets -G_{i1}$, which ensures that $x=0$
is feasible; third, $G_{i1}\gets -G_{i1} - 1$ to increase the
constraint distance from $x=0$.

Figure \ref{figRandLPs} compares the novel methods AugLag and
AnyAugLag with standard LogBarrier and SqrPenalty. All methods
reliably find the same optimum with very small constraint violation
$\sum_i [g_i(x)]_+$. Interestingly, the any-time augmented lagrangian
methods performs extremely well for moderate problem sizes, but
clearly looses its benefits for larger sizes. We inspected its
behavior and found qualitatively that the declined performance
coincides with significant non-stationarity of the constraint activity
also in the later stage of the random LP optimization. As anticipated
by our discussion and motivation of the proposed method, for random
LPs we should not expect stationarity of constraint activity during
optimization---infact, finding the set of active constraints is the
main problem for LPs and if we knew this set early the remaining
optimization would be trivial. The severe non-stationarity of
constraints works against the implicit assumptions made in the AnyAula
update, explaining its declining performance for many constraints.

\subsection{Robotic Trajectory Optimization}

\begin{figure}\centering
\label{figRobotEx}
\showh[.35]{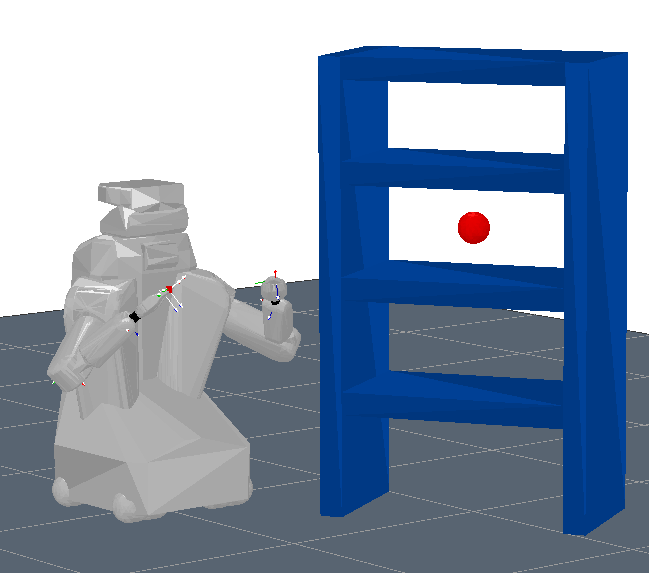}
\qquad
\showh[.3]{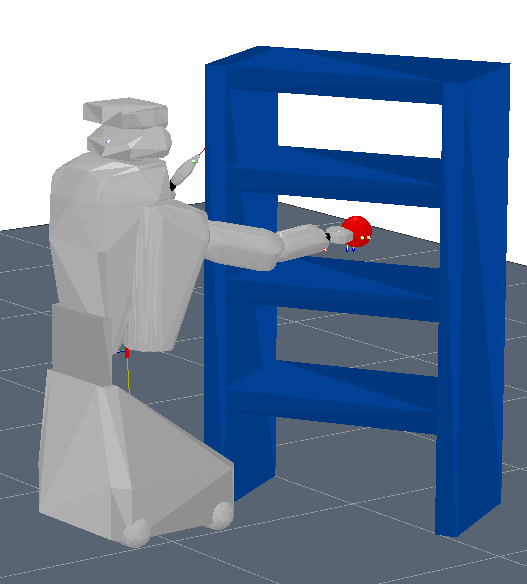}
\caption{Start and (optimized) goal configuration of a typical trajectory
optimization problem. The configuration space is $25$-dimensional,
the trajectory composed of 200 time slices, making this a
5000-dimensional problem over $x\in\RRR^{25\times200}$. The
computational cost is fully dominated by the number of evaluations of
$f(x)$ and $g(x)$, which implies computing potential collisions. }
\end{figure}

We also tested the performance on standard robotic trajectory
optimization problems, as illustrated in
Figure \ref{figRobotEx}. Before discussing the results, we would like
to characterize such problems and the methods we seek for: Feasible
path finding, let alone finding globally optimal paths, are in general
hard computational problems (NP-complete when discretizing the
configuration space). Therefore the approaches can roughly
be separated in two categories: path finding methods that are globally
(probabilistically) complete and local trajectory optimization methods
that aim to converge robustly and fast to a local feasible
optimum. Note that locality here is meant in the trajectory space, not
configuration space. Therefore, depending on the specifics of the cost
function and whether the optimization method allows to temporarily
traverse infeasible regions, local optimization can very well solve
problems that are in other contexts (potential fields) considered as
local deadlocks. In this view, here we aim for optimization methods
that robustly and agressively move towards local optima and, in the
vicinity of such local optima, efficiently minimize the local
non-linear convex problem.

Table \ref{tabRobotEx} displays the performance of the various
methods. Both, AnyAula and Aula converge to the same optimum (modulo
stopping criterion tolerance), whereas LogBarrier fails to find any
reasonable solution. In the light of the above discussion this can be
explained as follows: Finding a fully feasible path from the initial
trajectory $x$ to the final optimal trajectory $x^*$, where none of
the intermediate trajectories violates constraints, is very hard.
Aula and AnyAula implicitly relax the constraints in early iterations,
leading to much better convergence to a local optimum.

\begin{table}\centering
\label{tabRobotEx}
\begin{tabular}{|p{.3\columnwidth}||c|c|c|}
\hline
method
& $\l$ or $\mu$-updates
& $f$ evaluations
& suboptimality \\
\hline
AnyAula
& 20.25$\pm$2.3 & 48.25$\pm$4.93 & 0.05$\pm$0.03 \\
Aula
& 22.8$\pm$1.3 & 64.2$\pm$2.03 & 0.14$\pm$0.12 \\
LogBarrier
& 11$\pm$0 & 60.2$\pm$4.0 & 72337$\pm$3325 \\
SqrP
&11$\pm$0 & 42.6$\pm$3.2 & 4.45$\pm$0.91 \\
\hline
\end{tabular}
\caption{Performance on the robot trajectory optimization
problem. Averages are taken over randomization of the initial
configuration. Suboptimality denotes the difference in the found
minimum $f(x^*)$ to the best found by all methods.}
\end{table}

\small
\bibliography{refs,bibs}

\small
\bibliography{bibs}
\end{document}